\documentclass[11pt]{amsart}
\usepackage{epsfig}
\usepackage{graphicx}
\usepackage{color}

\newtheorem{theorem}{Theorem}[section]
\newtheorem{lemma}[theorem]{Lemma}

\newtheorem{cor}[theorem]{Corollary}

\newcommand{\lcr}{\raisebox{-5pt}{\mbox{}\hspace{1pt}
                  \epsfig{file=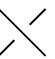}\hspace{1pt}\mbox{}}}
\newcommand{\ift}{\raisebox{-5pt}{\mbox{}\hspace{1pt}
                  \epsfig{file=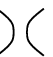}\hspace{1pt}\mbox{}}}
\newcommand{\zer}{\raisebox{-5pt}{\mbox{}\hspace{1pt}
                  \epsfig{file=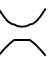}\hspace{1pt}\mbox{}}}
\theoremstyle{definition}

\theoremstyle{remark}

\numberwithin{equation}{section}

\title[The action of the Kauffman bracket skein algebra]{The action of the Kauffman bracket skein algebra of the torus on the
  Kauffman bracket skein module of the 3-twist knot complement}
\author{R{\u{a}}zvan Gelca}
\address{Department of Mathematics and Statistics, 
Texas Tech University, Lubbock, TX 79409}
\email{rgelca@gmail.com}
\author{Hongwei Wang}
\address{Department of Mathematics and Physics, Texas A\&M International University, Laredo, TX 78041}
\email{hongwei.wang@tamiu.edu}



\begin{document}
\maketitle

\begin{abstract}
We determine the action of the Kauffman bracket skein algebra of the torus
on the Kauffman bracket skein module of the complement of the 3-twist knot.
The point is to study the relationship between knot complements
and their boundary tori, an idea that has proved very fruitful in knot theory.
We place this idea in the context of Chern-Simons theory, where
 such  actions arose in connection with
the computation of the noncommutative version of the A-polynomial
that was defined in \cite{frgelo}, but they can also be interpreted as  quantum mechanical
systems. Our goal is to exhibit a detailed example in a part 
of Chern-Simons theory where examples are scarce. 
\end{abstract}

\section{Introduction}

This paper should be viewed as a piece of experimental mathematics.
It  describes the action of the Kauffman bracket skein algebra
of the torus on the Kauffman bracket skein module of the complement
of the 3-twist knot, which is listed as the $5_2$ knot in the knot table.
Such computations have been done before 
for the trefoil knot \cite{gelcatrefoil}, the figure-eight knot 
\cite{gelcasain}, and $(2,2p+1)$ torus knots \cite{nagasato} as the main step in the computation of the
noncommutative version of the A-polynomial  defined
in \cite{frgelo}. The noncommutative version of the A-polynomial
has been linked to colored Jones polynomials in \cite{frgelo}, \cite{gelcaproc},  \cite{garle}, and to  $SL(2,{\mathbb C})$-Chern-Simons theory \cite{gukov}, \cite{dimofte}, a
difficult area of mathematics that has yet to be thoroughly understood. 
The theory of Kauffman bracket skein modules has been linked to
$SL(2,{\mathbb C})$-character varieties \cite{bullock},
\cite{przytyckisikora}, as deformations of  rings of affine characters,
and as such they are also supposed to be related to $SL(2,{\mathbb
  C})$-Chern-Simons theory, though it is not known how. 

The case of the 3-twist knot  is probably the most complicated example
that can still be done by hand; this is why we want to present it to
the public. Computational complexity grows very fast
in the theory of skein modules, hence there are  few examples.
A striking feature  exhibited in this paper is the
ocurrence of Jones-Wenzl idempotents in the computations within skein
modules. This feature has been observed before by the first author;
it seems that (arbitrary)
skein computations tend to structure themselves in terms
of Jones-Wenzl idempotents.

In \cite{witten}, E. Witten has related $SU(2)$-Chern-Simons theory to
quantizations of moduli spaces of flat connections on surfaces,
which then leads to quantum mechanical models (see \cite{gelcauribe} for
a complete discussion). These quantum mechanical models have
a combinatorial version that arises from quantizing Wilson lines,
formulated using {\em reduced} skein modules (which are the building blocks
of the topological quantum field theory of Blanchet, Habbegger, Masbaum
and Vogel \cite{bhmv}). 
The action of the Kauffman bracket skein algebra of the torus
on the skein module of the knot complement that makes the object
of this paper is similar to that quantum mechanical model, but
here we work with the non-reduced version of skein modules. 
One might ask what does the model in which we work with the
actual skein modules and not their reduced versions represent?
Given such questions, the lack of  examples, and recent renewed interest
in skein modules,  we consider worth showing this particular situation,
as it might help clarify the general situation.

Given that the structure of the Kauffman bracket skein modules 
is now known for a fairly large family of knot and link complements
\cite{lepaper}, \cite{letran}, we hope that the above work will
be expanded to the study of structures that arise from skein
modules in other knot complements. 

We start with some background material.
Throughout the paper $t$  is a variable.
A {\em framed link} in an orientable 3-manifold $M$ is a 
disjoint union of embedded annuli. If $M$ is
the cylinder over the torus, framed links are identified with 
curves on the torus, with the annulus being parallel to the torus.
If we draw a framed link on paper, its framing is parallel to the plane of
the paper, unless the link is drawn on the torus, when we use the previous
convention. 
Let $\mathcal{L}$ be the set of isotopy classes of framed links in 
the manifold $M$, including the empty link. 
Consider the free ${\mathbb C}[t,t^{-1}]$-module with basis $\mathcal{L}$,
and factor it by the smallest subspace containing all expressions
of the form $\displaystyle{\lcr-t\zer-t^{-1}\ift}$
and 
$\bigcirc+t^2+t^{-2}$, where the links in each expression are
identical except in a ball in which they look like depicted.
This  quotient is denoted by $K_t(M)$ and is called the Kauffman 
bracket skein module of the manifold \cite{przytycki}. 
The factorization allows us to smoothen crossings (which we
can create at will using isotopy) and
to replace  trivial link components by a scalar. 
Because at each application of the  first skein relation
one term is replaced by two terms, the complexity of computations
grows exponentially, and so the computations in this paper are quite involved. 

For the cylinder over
torus, ${\mathbb T}^2\times I$ (where $I=[0,1]$), the skein module has a 
multiplication induced
by the operation of gluing one cylinder on top of another.
This multiplication has been explicated in \cite{frohmangelca},
here it is: As a module, $K_t({\mathbb T}^2\times I)$
is free with basis $(p,q)_T$, $p,q\in {\mathbb Z}$, $p\geq 0$,
where $(p,q)_T=T_n((p,q))$ with $T_n$ the (normalized)
Chebyshev polynomial of first kind defined
by  $T_0(x)=2$, $T_1(x)=x$, and $T_{n+1}(x)=xT_n(x)-T_{n-1}(x)$,
$n=\mbox{gcd}(p,q)$,
and $(p,q)$ is the curve of slope $ p/q$ on the torus. 
We have   product-to-sum formula
\begin{eqnarray*}
(p,q)_T*(r,s)_T=t^{|^{pq}_{rs}|}(p+r,q+s)_T+t^{-|^{pq}_{rs}|}(p-r,q-s)_T.
\end{eqnarray*}

In this paper we focus on the $3$-twist knot $K$ drawn in bold
line  in
Figure~\ref{twistknot}.  
\begin{figure}[h]
  \centering
\scalebox{.3}{\input{twistknot.pstex_t}}

\caption{}
\label{twistknot}
\end{figure}
When talking about the complement of $K$ we mean the compact orientable
manifold $S^3\backslash N(K)$ obtained by removing from the 3-sphere
an open regular neighborhood $N(K)$ of $K$.
The operation of gluing the cylinder over $\partial N(K)={\mathbb T}^2$
 to $S^3\backslash N(K)$ induces
a $K_t({\mathbb T}^2\times I)$-left module structure  on $K_t(S^3\backslash N(K))$.
In what follows we explicate this module structure.

It was shown in \cite{bullocklofaro} that
$K_t(S^3\backslash N(K))$ is a free ${\mathbb C}[t,t^{-1}]$-module
with basis $x^ny^k$,  $n\geq 0$, $0\leq k\leq 3$, where $x,y$
are shown in Figure~\ref{twistknot}. It suffices to
understand
the action of a set of generators of  $K_t({\mathbb T}^2\times I)$
on the basis,
and as generators we have chosen 
$(0,1)_T$, $(1,-3)_T$ and $(1,-2)_T$. The action of $(0,1)_T$ is 
$(0,1)_Tx^ny^k=x^{n+1}y^k$, so we  focus on the other two. 

Using the fact the $x$ can be pulled back into the cylinder over the boundary
as the skein $(0,1)_T$, and using the relations $(1,q)_T(0,1)_T=t(1,q+1)_T+t^{-1}(1,q-1)_T$,
and $(0,1)_T(1,q)_T=t^{-1}(1,q+1)_T+t(1,q-1)_T$,
we see that the action of $(1,q)_T$ on $x^ny^k$ can be found easily if we
know how $(1,q)_T$ acts on the basis elements $1=y^0,y,y^2,y^3$. It should
also be noted that in computations from this paper
$x$ behaves like a {\em scalar}. 

We will  change the basis of $K_t(S^3\backslash N(K))$ to
$S_n(x)S_k(y)$, where $S_n$ is the (normalized) Chebyshev polynomial of
the second kind: $S_0(x)=1, S_1(x)= x, S_{n+1}(x)=xS_n(x)-S_{n-1}(x)$. As such,
the basis elements are the curves $x$ and $y$ colored by {\em Jones-Wenzl idempotents}. There are two explanations for this, one is practical: the formulas
become simpler. But there is a deeper explanation for this, namely that
the polynomial $S_n$ is the character of the $n+1$-dimensional
irreducible representation of $SL(2,{\mathbb C})$, and as such the skein
$S_n(x)S_k(y)$
consists of two Wilson lines (one for $x$ and one for $y$) associated
to irreducible representations. It is worth pointing out that the
colored Jones polynomials of a knot $K$ are $(-1)^n\left<S_n(K)\right>$,
where $\left<\cdot\right>$ denotes the Kauffman bracket (of knots and links
in $S^3$).

In short, the goal of the paper is
to find $(1,-3)_T\cdot S_k(y)$ and $(1,-2)_T\cdot S_k(y)$,
 $k=0,1,2,3$.

\section{Formulas in a quotient of the Kauffman bracket skein module of cylinder over
  the twice punctured disk}\label{section2}

It is known that the Kauffman bracket skein module of the cylinder over the
twice punctured disk, i.e. a disk with  two disjoint
open disks being removed, is free with basis $x^my^nz^k$, $m,n,k\geq 0$,
 where $x$ and $z$ are curves that are parallel with the boundaries of
the two open disks that  have been removed, and $y$ is a curve parallel
to the boundary of the original disk. In Sections~\ref{section2} and \ref{section3}, we make the following convention. We schematically represent the
cylinder over the twice punctured disk sideways,
by drawing only the two curves that
trace the punctures in the cylinder. These curves will either be represented
as twisting around each other, such as in the first diagram from
Figure~\ref{skeinsXY}, or as two parallel lines such as in the second
and third diagram from the same figure. Closed curves  in the diagram
comprise skeins, taken with the blackboard framing. Whenever a number
is written next to a curve, such as the $k$ written next to the $y$-curve
in the first diagram from Figure~\ref{skeinsXY}, that number
indicates that the skein contains that many parallel copies of that curve,
as such as in our example there are $k$ parallel copies of $y$.


We factor the Kauffman bracket skein module of the cylinder over the twice punctured
disk by  the relation $$x=z$$ and perform all computations from this section
of the paper in this quotient. All computations in this section can
be used for general twist knots.

In Figure~\ref{skeinsXY}, we recall the skeins $X_i*y^k$ from \cite{gelcanagasato}
and  define the skeins $Y_1*y^k$. In the first diagram, the index $i$ counts
the crossings of the two strands that define the genus 2 handlebody. For $Y_1*y^k$, the
undercrossings can be at the bottom and the overcrossings 
at the top, as one diagram is mapped into the other by isotopy. 
\begin{figure}[h]
  \centering
\scalebox{.4}{\input{skeinsXY.pstex_t}}

\caption{}
\label{skeinsXY}
\end{figure}

\begin{lemma}\label{lemmaA} 
	The skeins $X_1*y^k$ and $Y_1*y^k$, $k\geq 0$, satisfy the recursions:
	\begin{eqnarray*}
		&&X_1*y^{k+1}=t^4yX_1*y^k+(t^{-2}-t^6)Y_1*y^k+2(1-t^4)x^2y^k,\quad k\geq 0,\\
		&&Y_1*y^{k+1}=t^{-4}yY_1*y^k+(t^2-t^{-6})X_1*y^k+2(1-t^{-4})x^2y^k,\quad k\geq 0,\\
		&&X_1*y^0=X_1=-t^4y-t^2x^2,\quad Y_1*y^0=Y_1=-t^2-t^{-2}.
	\end{eqnarray*}
	\end{lemma}

\begin{proof}
We start computing $X_1*y^{k+1}$ as in Figure~\ref{lemmaA1}.
\begin{figure}[h]
  \centering
\scalebox{.9}{\input{lemmaA1.pstex_t}}

\caption{}
\label{lemmaA1}
\end{figure}
The first diagram is computed as in Figure~\ref{lemmaA2new}, and is equal
to $t^2yX_1*y^k-t^4Y_1*y^k-2t^2x^2y^k$. So the first term
is $t^4yX_1*y^k-t^6Y_1*y^k-2t^4x^2y^k$. The sum of the  other 
terms is equal to $x^2y^k+x^2y^k+t^{-2}Y_1*y^k$. Adding we obtain the 
first recursion. Similarly for the second recursion. 
\begin{figure}[h]
  \centering
\scalebox{.7}{\input{lemmaA2new.pstex_t}}

\caption{}
\label{lemmaA2new}
\end{figure}
\end{proof}

Applying Lemma~\ref{lemmaA} we obtain
\begin{eqnarray*}
&&X_1*y=-t^8S_2(y)-t^6x^2S_1(y)+2(1-t^4)x^2+(t^4-1-t^{-4});\\
&&Y_1*y=-(t^6+t^{-6})S_1(y)+(2-t^4-t^{-4})x^2;
\end{eqnarray*}
\begin{eqnarray*}
&&X_1*y^2=-t^{12}S_3(y)+(-t^{10})x^2S_2(y)+(-2t^8+2)x^2S_1(y)\\
&&\quad +(t^8-2t^4-t^{-8})S_1(y)+(-2t^{6}+2t^{-2}-t^{-6})x^2;\\
&&Y_1*y^2=(-t^{10}-t^{-10})S_2(y)+(2-t^8-t^{-8})x^2S_1(y)\\
&&\quad+(-2t^6-2t^{-6}+2t^2+2t^{-2})x^2+(t^6+t^{-6}-2t^2-2t^{-2});
\end{eqnarray*}
The following result is a straighforward generalization of Lemma 1 in \cite{gelcanagasato}.

\begin{lemma}\label{lemma1}
	The skeins $X_i*y^k$, 
$i,k\geq 0$, satisfy the recursive relation
	\begin{eqnarray*}
&&		X_{i+2}*y^k=t^2yX_{i+1}*y^k-t^4X_i*y^k-2t^2x^2y^k,\\
&& X_{2}*y^k=t^2yX_{1}*y^k-t^4Y_1*y^k-2t^2x^2y^k.
	\end{eqnarray*}
\end{lemma}

As a consequence, we obtain
\begin{eqnarray*}
&&X_2=-t^6S_2(y)-t^4S_2(x)S_1(y)-t^4S_1(y)-2t^2S_2(x)-t^2;\\
&&X_3=-t^8S_3(y)-t^6S_2(x)S_2(y)-t^6S_2(y)-2t^4S_2(x)S_1(y)\\&&\quad
-t^4S_1(y)-2t^2S_2(x)-2t^2;\\
&&X_4=-t^{10}S_4(y)-t^8S_2(x)S_3(y)-t^8S_3(y)-2t^6S_2(x)S_2(y)\\
&&\quad -t^6S_2(y)-2t^4S_2(x)S_1(y) -2t^4S_1(y)-2t^2S_2(x)-2t^2;
\end{eqnarray*}
\begin{eqnarray*}
&&X_2*y=-t^{10}S_3(y)-t^8S_2(x)S_2(y)-t^8S_2(y)-2t^6S_2(x)S_1(y)\\
&&\quad +(-t^6-t^2)S_1(y)+(-2t^4+1)S_2(x)+(-2t^4+1);\\
  &&X_3*y=-t^{12}S_4(y)-t^{10}S_2(x)S_3(y)-t^{10}S_3(y)-2t^8S_2(x)S_2(y)\\
&&\quad +(-t^8-t^4)S_2(y)+(-2t^6-t^2)S_2(x)S_1(y)+(-2t^6-t^2)S_1(y)\\
&&\quad -2t^4S_2(x) +(-2t^4+1);\\
&&X_4*y=-t^{14}S_5(y)-t^{12}S_2(x)S_4(y)-t^{12}S_4(y)-2t^{10}S_2(x)S_3(y)\\
&&\quad +(-t^{10}-t^6)S_3(y)+(-2t^8-t^4)S_2(x)S_2(y)+(-2t^8-t^4)S_2(y)\\
&&\quad +(-2t^6-2t^2)S_2(x)S_1(y)+(-2t^6-t^2)S_1(y)-2t^4S_2(x)-2t^4;
\end{eqnarray*}
\begin{eqnarray*}
&&X_2*y^2=-t^{14}S_4(y)-t^{12}S_2(x)S_3(y)-t^{-12}S_2(y)-2t^{10}S_2(x)S_2(y)\\
&&\quad
  +(-t^{10}-2t^6)S_2(y)+(-2t^4-2t^8+2)S_2(x)S_1(y)\\
	&&\quad +(-2t^4-2t^8+2)S_1(y)+(-2t^6-2t^2+2t^{-2})S_2(x)\\
	&&\quad +(-2t^6+t^{-2}-t^{-6});\\
&&X_3*y^2=
  -t^{16}S_5(y)-t^{14}S_2(x)S_4(y)-t^{14}S_4(y)+(-2t^{12})S_2(x)S_3(y)\\
&&\quad +(-t^{12}-2t^8)S_3(y)+(-2t^{10}-2t^6)S_2(x)S_2(y)+(-2t^{10}-2t^6)S_2(y)\\
&&\quad +(-2t^8-4t^4+2)S_2(x)S_1(y)+(-2t^8-2t^4+1)S_1(y)\\
&&\quad+(-2t^6-2t^2+t^{-2})S_2(x)+(-2t^6-2t^2+t^{-2});\\
&&X_4*y^2=-t^{18}S_6(y)-t^{16}S_2(x)S_5(y)-t^{16}S_5(y)-2t^{14}S_2(x)S_4(y)\\
&&\quad +(-t^{14}-2t^{10})S_4(y)+(-2t^{12}-2t^8)S_2(x)S_3(y)\\
&&\quad +(-2t^{12}-2t^8)S_3(y)+(-2t^{10}-4t^6)S_2(x)S_2(y)\\
&&\quad +(-2t^{10}-2t^6-t^2)S_2(y)+(-2t^8-4t^4+1)S_2(x)S_1(y)\\
&&\quad +(-2t^8-4t^4+1)S_1(y)+(-2t^6-2t^2)S_2(x)+(-2t^6-2t^2+t^{-2}).
\end{eqnarray*}


\begin{lemma}\label{lemma2}
	The skeins $X_i*y^k$, $i,k\geq 0$, satisfy the recursive relation
	\begin{eqnarray*}
		X_2*y^k=t^{-2}X_1*y^{k+1}-2t^{-2}x^2y^k-t^{-4}Y_1*y^k
	\end{eqnarray*}
	and for $i\geq 1$, 
	\begin{eqnarray*}
		X_{i+2}*y^k=t^{-2}X_{i+1}*y^{k+1}-2t^{-2}x^2y^k-t^{-4}X_{i}*y^k.
	\end{eqnarray*}
\end{lemma}

\begin{proof}
To compute $X_i*y^{k+1}$, we separate a $y$ from $y^{k+1}$, slide it so as to produce two crossings
in the link diagram, then solve the crossings using the Kauffman
bracket skein relation as  in  Figure~\ref{lemma2figure}.  
\begin{figure}[h]
  \centering
\scalebox{.4}{\input{lemma2.pstex_t}}

\caption{}
\label{lemma2figure}
\end{figure}

In the first term, by sliding  the strand to the right we see that this
 term equals $t^2X_{i+1}*y^k$. The second and  third terms
are each equal to $2x^2y^2$. The last term is equal to 
$X_{i-1}*y^k$ if $i\geq 2$ and to $Y_1*y^k$ if $i=1$. 
\end{proof}

Define the skeins $A*y^k$, $\overline{A*y^k}$,  $B*y^k$,  $\overline{B*y^k}$ 
as in Figure~\ref{asonly}.
\begin{figure}[h]
  \centering
\scalebox{.35}{\input{asonly.pstex_t}}

\caption{}
\label{asonly}
\end{figure}

\begin{lemma}\label{recasbs}
The following relations hold
  \begin{eqnarray*}
  &&A*S_k(y)=(-t^{2k+2}-t^{-2k-2})B*S_k(y),\\
  &&\overline{A*S_k(y)}=(-t^{2k+2}-t^{-2k-2})\overline{B*S_k(y)}\\
  &&B*y^k=t^2yB*y^{k-1}+(1-t^{-4})\overline{B*y^{k-1}},\\
  &&\overline{B*y^k}=t^{-2}y\overline{B*y^{k-1}}+(1-t^4)B*y^{k-1},\\
&&B*y^0=\overline{B*y^0}=x.
\end{eqnarray*} 
\end{lemma}

\begin{proof}
  The formulas for $A*S_k(y)$, $\overline{A*S_k(y)}$ follow from the standard
  properties of Jones-Wenzl idempotents.



For  $B*y^k$ (see Figure~\ref{aj}) resolve
the crossings specified by the arrow to obtain the first sum in this figure.
Perform an isotopy of the first skein from the sum to obtain
the first skein on the second row (in the process we remove and then add
a positive twist), then remove a negative twist from the second term
and perform an isotopy in this term. Then  apply the 
Kauffman bracket skein relation in the place specified by the arrow
to obtain the desired relation.

  $\overline{B*y^k}$ is obtained by reflecting  $B*y^k$ over
a horizontal line, and under reflections, in the  Kauffman
bracket  $t$ is replaced by $t^{-1}$.  
\begin{figure}[h]
  \centering
\scalebox{.35}{\input{aj.pstex_t}}

\caption{}
\label{aj}
\end{figure}
\end{proof}

\begin{cor}
The following formulas hold
\begin{eqnarray*}
&&A*y^0=(-t^2-t^{-2})x,\\
&&A*y=(-t^6-t^{-2})xS_1(y)+(-t^4+1-t^{-4}+t^{-8})x,\\
&&A*y^2=(-t^{10}-t^{-2})xS_2(y)+(-t^8+1-t^{-4}+t^{-12})xS_1(y)\\
&&\quad +(-t^6-t^{-6}-t^{-2}+t^{-10})x\\
&&A*y^3=(-t^{14}-t^{-2})xS_3(y)+(-t^{12}+1-t^{-4}+t^{-16})xS_2(y)\\
&&\quad+(-t^{10}-2t^{-2}+t^2-t^{-6}+t^{-14}-2t^6)xS_1(y)\\&&\quad+(-t^8-t^4+2+t^{-12}-2t^{-4}+t^{-8})x.
\end{eqnarray*}
\end{cor}

\begin{cor}
The following formulas hold
\begin{eqnarray*}
&&\overline{A*y^0}=(-t^2-t^{-2})x,\\
&&\overline{A*y}=(-t^2-t^{-6})xy+(t^{8}-t^4+1-t^{-4})x,\\
&&\overline{A*y^2}=(-t^{2}-t^{-10})xS_2(y)+(t^{12}-t^{4}+1-t^{-8})xS_1(y)\\
&&\quad +
(-t^6-t^{-6}-t^{2}+t^{10})x\\
&&\overline{A*y^3}=(-t^{2}-t^{-14})xS_3(y)+(t^{16}-t^{4}+1-t^{-12})xS_2(y)\\
&&\quad+
(-t^{-10}+t^{-2}-2t^2-t^{6}+t^{14}-2t^{-6})xS_1(y)\\&&\quad
+(-t^{-8}-t^{-4}+2+t^{12}-2t^4+t^{8})x.\\
\end{eqnarray*}
\end{cor}

\begin{proof}
The skein $\overline{A_0*y^k}$ is the reflection of  $A_0*y^k$
over a horizontal line. To get the  formulas for
 $\overline{A_0*y^k}$, swap $t$ and $t^{-1}$ in the 
formulas for  $A_0*y^k$.
\end{proof}

 We define the skeins $C_j*y^k, D_j*y^k, E_j, F_j$ as in
Figure~\ref{cdef2}. 
\begin{figure}[h]
  \centering
\scalebox{.35}{\input{cdef2.pstex_t}}

\caption{}
\label{cdef2}
\end{figure}

\begin{lemma}\label{csdsesfs}
The skeins $C_j*y^k, D_j*y^k, E*y^j,F*y^j$ satisfy the following relations
\begin{eqnarray*}
&&C_j*y^k=t^2C_{j+1}*y^{k-1}+(1-t^{-4})D_j*y^{k-1}\\
&&D_j*y^k=t^{-2}D_{j+1}*y^{k-1}+(1-t^4)C_j*y^{k-1}\\
&&C_j*y^0=x\overline{B*y^j}, \quad D_j*y^0=t^{-2}x\overline{B*y^j}+(1-t^4)E*y^j\\
  &&E*S_j(y)=(-t^{2j+2}-t^{-2j-2})S_j(y),  F*S_j(y)=(t^{2j+2}+t^{-2j-2})^{2}S_j(y).
\end{eqnarray*} 
\end{lemma}

\begin{proof}
For $C_j*y^k$, separate
a strand from $y^k$ as in the skein on  the left in Figure~\ref{csds}, then
resolve the crossings defined by arrows to obtain 
\begin{eqnarray*}
t^2C_{j+1}*y^{k-1}+D_j*y^{k-1}+D_j*y^{k-1}+t^{-2}(-t^2-t^{-2})D_j*y^{k-1},
\end{eqnarray*}
which yields the relation.  For $D_j*y^k$ do  the same  in the skein on the right. 
\begin{figure}[h]
  \centering
\scalebox{.35}{\input{csds.pstex_t}}

\caption{}
\label{csds}
\end{figure}

The skein $C_j*y^0$ is the mirror image of $xA_j$ over a vertical line,
so, as a skein, equals $x\overline{A_j}$. Resolving the two crossings in $D_j*y^0$ specified in 
Figure~\ref{dande}, we get $t^2(-t^2-t^{-2})E_j+E_j+E_j+t^{-2}x\overline{A_j}$.
 \begin{figure}[h]
   \centering
\scalebox{.35}{\input{dande.pstex_t}}

\caption{}
\label{dande}
\end{figure}
 The formulas for $E*S_j(y)$ and $F*S_j(y)$ follow from standard properties of
 Jones-Wenzl idempotents.

\end{proof}

\begin{cor}
The following formulas hold
\begin{eqnarray*}
&&C_0*y^0=x^2\\
&&C_0*y=x^2S_1(y)+(t^{-2}+t^2-t^6-t^{-6})x^2+(t^6-t^{-2}-t^2+t^{-6})\\
&&C_0*y^2=x^2S_2(y)+(t^2+t^{-2}-t^{10}-t^{-10})x^2S_1(y)+(3-t^{8}-t^{-8})x^2\\
&&\quad +(-t^{-6}-t^6+t^{10}+t^{-10})S_1(y)\\
&&C_0*y^3=x^2S_3(y)+(t^{-2}+t^2-t^{14}-t^{-14})x^2S_2(y)\\
&&\quad +(t^4+t^{-4}+4-t^{12}-t^{-12}-t^{-8}-t^8)x^2S_1(y)\\&&\quad 
+(-2t^{-6}-2t^{10}+4t^{-2}+5t^2+2t^{14}-3t^6-t^{-10}+t^{-14})x^2\\
&&\quad +(-t^{10}+t^{14}-t^{-10}+t^{-14})S_2(y)\\&&\quad +(-3t^{-2}+3t^{-6}-3t^2+3t^6).
\end{eqnarray*}
\end{cor}

We define the skeins  $G*y^k,H*y^k$, and $\overline{H*y^k}$ as in
Figure~\ref{gshsjs}. The next result has a proof analogous to that
to Lemma~\ref{recasbs}.
\begin{figure}[h]
  \centering
\scalebox{.35}{\input{gshsjs.pstex_t}}

\caption{}
\label{gshsjs}
\end{figure}

\begin{lemma}\label{gshslemma}
The following formulas hold
\begin{eqnarray*}
&&G*S_k(y)=(-t^{2k+2}-t^{-2k-2})H*S_k(y),\\
  &&H*y^k=t^2yH*y^{k-1}+(1-t^{-4})J*y^{k-1}, \quad H*y^0=y,\\
  && J*S_k(y)=(-t^{2k+2}-t^{-2k-2})S_k(y).
\end{eqnarray*}
\end{lemma}




\begin{cor}
The following formulas hold
\begin{eqnarray*}
&&G*y^0=(-t^{2}-t^{-2})S_1(y),\\
  &&G*y=(-t^6-t^{-2})S_2(y)+(-t^{-2}-t^{-10}),\\
  &&G*y^2=(-t^{10}-t^{-2})S_3(y)+(-t^2-2t^{-2}-t^{-14})S_1(y),\\
  &&G*y^3=(-t^{14}-t^{-2})S_4(y)+(-2t^6-3t^{-2}-t^{-18})S_2(y)+(-2t^{-2}-2t^{-10}).
\end{eqnarray*}
\end{cor}

\section{Formulas in the Kauffman bracket skein module of the 
3-twist knot complement}\label{section3}

For the complement of the 3-twist knot
Lemma~6 in \cite{gelcanagasato} gives

\begin{lemma}\label{lemma6} For all $k\geq 0$ we have
\begin{eqnarray*}
X_4*y^k=-t^{-4}X_3*y^k-t^{-2}x^2y^k.
\end{eqnarray*}
\end{lemma}

Lemma~\ref{lemma6} yields different formulas for
$X_4*y^k$, which, combined with those in \S~2,
give relations that  successively compute $S_4(y)$, 
 $S_5(y)$, $S_6(y)$:
\begin{eqnarray*}
&&S_4(y)=[-t^{-2}S_2(x)-(t^{-2}+t^{-6})]S_3(y)+[-(2t^{-4}+t^{-8})S_2(x)\\
&&\quad -
(t^{-8}+t^{-4})]S_2(y) +[-2(t^{-10}+t^{-6})S_2(x)-(t^{-10}+2t^{-6})]S_1(y)\\
&&\quad+[-(t^{-12}+2t^{-8})S_2(x)-(t^{-12}+2t^{-8})];
\end{eqnarray*}
\begin{eqnarray*}
&&S_5(y)=[(t^{-4}S_4(x)+(t^{-4}+t^{-8})S_2(x)+(t^{-4}+t^{-12})]S_3(y)\\
&&\quad +[(2t^{-6}+t^{-10})S_4(x)+(3t^{-6}+2t^{-10}+t^{-14})
S_2(x)+(t^{-6}+t^{-10})]S_2(y)\\
&&\quad +[(2t^{-12}+2t^{-8})S_4(x)+(3t^{-12}+4t^{-8}+2t^{-16})S_2(x)\\
&&\quad+(2t^{-12}+2t^{-8}+t^{-16})]S_1(y)+(t^{-14}+2t^{-10})S_4(x)\\
&&\quad +(3t^{-14}+4t^{-10}+t^{-18})S_2(x)+(2t^{-18}+2t^{-14}+2t^{-10});
\end{eqnarray*}
\begin{eqnarray*}
&&S_6(y)=[-t^{-6}S_6(x)+(-t^{-6}-t^{-10})S_4(x)
+(-t^{-6}-t^{-14})S_2(x)\\
&&\quad +(-t^{-6}-t^{-18})]S_3(y)+[(-2t^{-8}-t^{-12})S_6(x)\\
&&\quad +(-3t^{-8}-2t^{-12}-t^{-16})
    S_4(x) +(-3t^{-8}-t^{-12}-t^{-16})S_2(x)+
\\
&&\quad +(-t^{-20}-t^{-16}-t^{-8})]S_2(y)+[(-2t^{-14}-2t^{-10})S_6(x)\\
&&\quad +(-4t^{-10}-3t^{-14}-2t^{-18})S_4(x)
    +(-2t^{-14}-4t^{-10}-2t^{-18})S_2(x)
\end{eqnarray*}
\begin{eqnarray*}
&&\quad +(-2t^{-10}-t^{-18})]S_1(y)+ [(-2t^{-12}-t^{-16})S_6(x)\\
&&\quad +(-4t^{-12}-3t^{-16}-t^{-20})S_4(x)
 +(-4t^{-12}-2t^{-16}-3t^{-20})S_2(x)\\
&&\quad +(-t^{-20}-2t^{-12}-t^{-24})].
\end{eqnarray*}

Also, from Lemma~\ref{lemma6}, we obtain
\begin{eqnarray*}
X_4=t^4S_3(y)+(t^2S_2(x)+t^2)S_2(y)+(2S_2(x)+1)S_1(y)+t^{-2}S_2(x)+t^{-2}.
\end{eqnarray*}

Combining  Lemma~\ref{lemma2} and Lemma~\ref{lemma6}
 we obtain the following recursive scheme
that allows the writing of $X_j*S_k(y)$, $1\leq j\leq 4,1\leq k\leq 6$ in terms
 of the basis $S_j(x)S_k(y)$, $0\leq j, 0\leq k\leq 3$: \\
\scalebox{.92}{$X_1*S_{k+1}(y)=t^2X_2*S_k(y)+t^{-2}Y_1*S_k(y)-X_1*S_{k-1}(y)+(2S_2(x)+2)S_k(y),$}\\
\scalebox{.92}{$X_2*S_{k+1}(y)=t^2X_3*S_k(y)+t^{-2}X_1*S_k(y)-X_2*S_{k-1}(y)+(2S_2(x)+2)S_k(y),$}\\
\scalebox{.92}{$X_3*S_{k+1}(y)=t^2X_4*S_k(y)+t^{-2}X_2*S_k(y)-X_3*S_{k-1}(y)+(2S_2(x)+2)S_k(y),$}\\
\scalebox{.92}{$X_4*S_k(y)=-t^{-4}X_3*S_k(y)-t^2x^2S_k(y).$}\\
Using also the formulas for $X_1,X_2,X_3,X_4$ and those that express $S_4(y)$,
$S_5(y)$, and $S_6(y)$ in terms of the basis
we obtain 
\begin{eqnarray*}
&&X_3*S_1(y)=t^6S_3(y)+t^4S_2(x)S_2(y)+t^2S_2(x)S_1(y)+S_2(x)+2;\\
&&X_4*S_1(y)=-t^2S_3(y)-S_2(x)S_2(y)+(-2t^{-2}S_2(x)-t^{-2})S_1(y)\\
&&\quad -t^{-4}S_2(x)-2t^{-4}.
\end{eqnarray*}
\begin{eqnarray*}
&&X_2*S_2(y)=t^8S_3(y)+t^6S_2(x)S_2(y)+[(2+t^4)S_2(x)+2]S_1(y)\\
&&\quad +(t^2+2t^{-2})S_2(x)+(2t^2+t^{-2}-t^{-6});\\
&&X_3*S_2(y)=-t^{4}S_3(y)-t^2S_2(x)S_2(y)-t^{-2};\\
&&X_4*S_2(y)=S_3(y)-t^{-2}S_2(y)+t^{-6}.
\end{eqnarray*}
\begin{eqnarray*}
&&X_2*S_3(y)=-t^6S_3(y)+(-t^4+1)S_2(x)S_2(y)+2S_2(y)+[2t^{-2}S_2(x)\\
&&\quad +(-2t^{-2}-t^{-10})]S_1(y)+(2t^{-4}-t^{-8})S_2(x)+(-1+2t^{-4}-t^8);\\
&&X_3*S_3(y)=t^2S_3(y)+2S_2(x)S_2(y)+S_2(y)+[2t^{-2}S_2(x)+2t^{-2}]S_1(y);\\
&&X_4*S_3(y)=[-t^{-2}S_2(x)-2t^{-2}]S_3(y)+[-2t^{-4}S_2(x)-t^{-4}]S_2(y)\\
&&\quad +[-2t^{-6}S_2(x)-2t^{-6}]S_1(y);
\end{eqnarray*}
\begin{eqnarray*}
&&X_2*S_4(y)=[2S_2(x)+(t^4+2)]S_3(y)+[(2t^2+2t^{-2})S_2(x)\\
&&\quad+(t^2-t^{-14})]S_2(y)+[(2+2t^{-4}-t^{-12})S_2(x)\\
&&\quad +(2+2t^{-4}-t^{-12})]S_1(y)+(2+2t^{-6}-2t^{-10})S_2(x)\\
&&\quad +(2t^{-2}+t^{-6}-t^{-10});\\
&&X_3*S_4(y)=S_2(x)S_3(y)+t^{-2}S_2(y)-t^{-12}S_1(y)-t^{-10}S_2(x);\\
&&X_4*S_4(y)=[t^{-4}S_4(x)+(2t^{-4}+t^{-8})S_2(x)+(2t^{-4}+t^{-8})]S_3(y)\\
&&\quad
  +[(2t^{-6}+2t^{-10})S_4(x)+(5t^{-6}+6t^{-10})S_2(x)+(2t^{-6}+4t^{-10})]S_2(y)\\
&&\quad+[(2t^{-8}+2t^{-12})S_4(x)+(10t^{-8}+5t^{-12})S_2(x)+(3t^{-12}+t^{-16})]S_1(y)\\
&&\quad +(2t^{-10}+t^{-14})S_4(x)+(6t^{-10}+3t^{-14})S_2(x)+(4t^{-10}+2t^{-14}).
\end{eqnarray*}

\section{The action of the skeins $(1,-3)_T$ and $(1,-2)_T$
 on the skein module of the 3-twist knot complement} 

As said in the introduction, we compute the action of 
 $(1,-3)_T$, $(1,-2)_T$ from $K_t({\mathbb T}^2\times I)$ 
on the basis elements $S_k(y)$, $0\leq k\leq 3$ of $K_t(S^3\backslash N(K))$.
The skeins $(1,-3)_T$ and $(1,-2)_T$ are depicted in Figure~\ref{twocurves},
with the cylinder ${\mathbb T}^2\times [0,1]$ embedded as a regular
neighborhood
of the boundary of the knot complement. Before starting the computation we prove a lemma.

\begin{figure}[h]
  \centering
\scalebox{.35}{\input{twocurves.pstex_t}}

\caption{}
\label{twocurves}
\end{figure}

\begin{lemma}\label{strand}
The identities from Figure~\ref{strandfig} hold.
Here the curved strand can encircle several parallel straight
strands. 
\end{lemma}

\begin{figure}[h]
  \centering
\scalebox{.35}{\input{strandfig.pstex_t}}

\caption{}
\label{strandfig}
\end{figure}

\begin{proof}
Pull the strand in the term on the left until you
create a negative twist in the first identity and a positive
twist in the second identity, then resolve the crossing. 
\end{proof}

First we find  $(1,-3)_T\cdot y^k$, $k=0,1,2,3$. For that we add  $y^k$ to the skein represented by
the curve on the left side of Figure~\ref{twocurves}, then push this curve inside the 
knot complement. There is a small technical detail. The framing that the curve inherits
from the torus does not coincide with the framing defined by the plane of the paper. The resulting skein
 (with framing defined by the plane of the paper) is  the skein  
 from Figure~\ref{firstcurve} multiplied by $t^6$. We  compute this 
skein from the figure first,  then multiply by the adjusting factor
in the end. In Figure~\ref{firstcurve}  we have labeled the 5 crossings
 in the order in which they are resolved. We use a boldface
curve for $y^k$, and here and in subsequent figures no longer write
the label $y^k$ next to it. 

\begin{figure}[h]
  \centering
\scalebox{.4}{\input{firstcurve.pstex_t}}

\caption{}
\label{firstcurve}
\end{figure}

We denote by a string of length $k$ consisting
of 
$+$'s and $-$'s inside double brackets the skein obtained from
$(1,-3)_T\cdot y^k$ by smoothening the first $k$ crossings (in the order of labels),
horizontally for a plus and vertically for a minus.
For example the Kauffman bracket skein relation applied to the first
crossing reads $(1,-3)_T  \cdot y^k= t((+))+t^{-1}((-))$.
Applying the Kauffman bracket skein relation repeatedly we obtain 
\begin{eqnarray*}
(1,-3)_T\cdot y^k=t^4((++++))+t^3((+++-+))+t((+++--))\\
+t((++-))+((+-))+t^{-1}((-))
\end{eqnarray*} 
where the skeins from this expression
are shown in Figure~\ref{plusminus}. 

\begin{figure}[h]
  \centering

\scalebox{.25}{\input{plusminus1.pstex_t}}
\scalebox{.25}{\input{plusminus2.pstex_t}}
\scalebox{.25}{\input{plusminus3.pstex_t}}

\scalebox{.25}{\input{plusminus4.pstex_t}}
\scalebox{.25}{\input{plusminus5.pstex_t}}
\scalebox{.25}{\input{plusminus6.pstex_t}}

\caption{}
\label{plusminus}
\end{figure}

After removing two negative twists in $((-))$ and focusing on the lower
part of the twist knot only,  we obtain that
$((-))$ is equal to the first skein in Figure~\ref{minus} multiplied
by $t^{-6}$. This new skein is computed as shown in the figure
by using the skein relation and sliding the strands.  In this
sum the first term is $t^2xy^{k+1}$, the
third is $xy^k$, and the fourth is $t^{-2}xX_3*y^k$. 
\begin{figure}[h]
  \centering
\scalebox{.25}{\input{minus.pstex_t}}

\caption{}
\label{minus}
\end{figure}

Let us focus on the second term in the sum. Applying  Lemma~\ref{strand}
in the place specified by the arrow, we can transform
it as in Figure~\ref{min2}. The first term is just $-t^2xy^{k+1}$.
By applying Lemma~\ref{strand} at the two places specified by
arrows we can transform the first  term into
$-t^4x^3y^k-t^4xy^k-t^6xX_4*y^k-t^2xX_3*y^k$. So the term
 we are computing equals
$-t^2xy^{k+1}-t^4x^3y^k-t^4y^k-t^6xX_4*y^k-t^2xX_3*y^k$. 
 
\begin{figure}[h]
  \centering
\scalebox{.25}{\input{min2.pstex_t}}

\caption{}
\label{min2}
\end{figure}

Therefore
$((-))=-xX_4*y^{k}+(t^{-8}-t^{-4})xX_3*y^{k}-t^{-2}x^3y^k+(t^{-6}-t^{-2})xy^k$, which after applying Lemma~\ref{lemma6} becomes
\begin{eqnarray*}
((-))=t^{-8}xX_3*y^k+(t^{-6}-t^{-2})xy^k.
\end{eqnarray*}

To compute $((++-))$ we look again at the lower part of
the twist knot and apply Lemma~\ref{strand} in the places 
specified by arrows in Figure~\ref{ppm}. In the last sum, by
resolving the two crossings
in each diagram 
by the skein relation we find that the first term is
\begin{eqnarray*}
&&t^{-4}xX_3*y^{k+1}+t^{-2}x^3y^{k+1}+t^{-4}xX_3*y^{k+1}+(-t^{-4}-t^{-8})xX_3*y^{k+1}\\
&&\quad=(t^{-4}-t^{-8})xX_3*y^{k+1}+t^{-2}x^3y^{k+1}.  
\end{eqnarray*}
\begin{figure}[h]
  \centering
\scalebox{.25}{\input{ppm.pstex_t}}

\caption{}
\label{ppm}
\end{figure}

Resolving the crossings in the second term we obtain
that it is equal to 
\begin{eqnarray*}
t^{-4}xy^{k+1}X_1+t^{-6}A+t^{-6}B+t^{-8}xX_3*y^{k+1},
\end{eqnarray*}
where skeins $A$ and $B$ are as in Figure~\ref{ppm2}.
\begin{figure}[h]
  \centering
\scalebox{.3}{\input{ppm2.pstex_t}}

\caption{}
\label{ppm2}
\end{figure}
Compute $B$ by applying 
Lemma~\ref{strand} at the location specified by  arrow to obtain $B=-t^2xX_3*y^{k+1}-t^4xy^{k+1}$.
Then transform $A$ by an isotopy,  use Lemma~\ref{strand} as shown in
Figure~\ref{ppm3} to obtain $A=t^{-2}xX_2*y^{k+1}+t^2xy^{k+2}+x^3y^{k+1}+xy^{k+1}$. 
Substitute $X_1$ by $-t^4y-t^2x^2$ to
conclude that the second term in the three-term sum from the second line
in  Figure~\ref{ppm} equals
\begin{eqnarray*}
(t^{-8}-t^{-4})xX_3*y^{k+1}+t^{-8}xX_2*y^{k+1}+(t^{-4}-1)xy^{k+2}\\
+(t^{-6}-t^{-2})x^3y^{k+1}+(t^{-6}-t^{-2})xy^{k+1}.
\end{eqnarray*}
 Similarly, the third term is
$(1-t^{-4})xy^k-t^{-6}xX_3*y^{k}$. Hence
\begin{eqnarray*}
&&((++-))=-t^{-6}xX_3*y^k+t^{-8}xX_2*y^{k+1}+(t^{-4}-1)xy^{k+2}\\
&&\quad +t^{-6}x^3y^{k+1}+(t^{-6}-t^{-2})xy^{k+1}+(1-t^{-4})xy^k.
\end{eqnarray*}
After applying Lemma~\ref{lemma2} this becomes
\begin{eqnarray*}
  &&((++-))=t^{-10}xX_1*y^k+(t^{-4}-1)xy^{k+2}+t^{-6}x^3y^{k+1}\\
 && \quad +(t^{-6}-t^{-2})xy^{k+1} +2t^{-8}x^3y^k+(1-t^{-4})xy^k.
  \end{eqnarray*}
\begin{figure}[h]
  \centering
\scalebox{.25}{\input{ppm3.pstex_t}}

\caption{}
\label{ppm3}
\end{figure}
We turn to   $((+-))$ and by working in  the lower
part of the knot,  we apply Lemma~\ref{strand} as
specified by an arrow in Figure~\ref{pm}, then remove
a negative twist in each  resulting diagram,
to obtain the two-term sum in Figure~\ref{pm}. 
A close examination shows that the last diagram in the figure
is the same as the last diagram in Figure~\ref{ppm}. Adjusting
for the different coefficient, we deduce that the second term from  the sum in Figure~\ref{pm} is 
$(-t^{-3}+t^{-7})xy^k+t^{-9}xX_3*y^{k}$. 
Computing similarly by resolving both crossings with
the Kauffman bracket skein relation, we find that the
first term is $t^{-7}xX_3*y^{k+1}+(t^{-5}-t^{-1})xy^{k+1}$. 
Combining, we get
\begin{eqnarray*}
((+-))=t^{-7}xX_3*y^{k+1}+t^{-9}xX_3*y^k+(t^{-5}-t^{-1})xy^{k+1}+(t^{-7}-t^{-3})xy^k.
\end{eqnarray*}
\begin{figure}[h]
  \centering
\scalebox{.25}{\input{pm.pstex_t}}

\caption{}
\label{pm}
\end{figure}

Now  turn to the computation of $((+++--))$.
After  an isotopy at the top part of the twist knot
we obtain the first diagram in Figure~\ref{pppmm1}, where in
this diagram we  ignore the top part of the twist knot, as we will not
use it again.  We apply Lemma~\ref{strand}
as specified by the arrow to obtain the sum on the right.
\begin{figure}[h]
  \centering
\scalebox{.25}{\input{pppmm1.pstex_t}}

\caption{}
\label{pppmm1}
\end{figure}
The first term is computed by applying Lemma~\ref{strand} at the point
specified by the arrow, as  in Figure~\ref{pppmm2}. 
In the sum from Figure~\ref{pppmm2}, the second term is 
$t^{-6}xY_1*y^k$. For the first term, we
perform an isotopy to make it look like in Figure~\ref{pppmm3},
then apply Lemma~\ref{strand} to obtain that it
is equal to $-t^{-6}x^3y^{k+1}-t^{-8}xyX_1*y^k$. 
\begin{figure}[h]
  \centering
\scalebox{.25}{\input{pppmm2.pstex_t}}

\caption{}
\label{pppmm2}
\end{figure}
Thus the 
first term of the sum in Figure~\ref{pppmm1} equals
 $-t^{-6}x^3y^{k+1}-t^{-8}xyX_1*y^k+t^{-6}xY_1*y^k$. 

\begin{figure}[h]
  \centering
\scalebox{.25}{\input{pppmm3.pstex_t}}

\caption{}
\label{pppmm3}
\end{figure}

The second term in Figure~\ref{pppmm1} can be transformed by an isotopy
into the first skein in Figure~\ref{pppmm5}. Then  apply Lemma~\ref{strand}
as specified by the arrow to obtain the sum on the right. The second term
is $t^{-8}xy^k$.  The first term can be computed by applying 
the Lemma~\ref{strand} as specified, and is 
$-t^{-8}x^3y^k-t^{-10}X_1*y^k. $

\begin{figure}[h]
  \centering
\scalebox{.25}{\input{pppmm5.pstex_t}}

\caption{}
\label{pppmm5}
\end{figure}

Combining the results we obtain
\begin{eqnarray*}
((+++--))=-t^{-6}x^3y^{k+1}-t^{-8}xyX_1*y^k+t^{-6}xY_1*y^k-t^{-8}x^3y^k\\
-t^{-10}xX_1*y^k+t^{-8}xy^k.
\end{eqnarray*}
Using Lemma~\ref{lemma1} we write this as
\begin{eqnarray*}
  ((+++--))=-t^{-10}xX_2*y^k-t^{-10}xX_1*y^k-t^{-6}x^3y^{k+1}-3t^{-8}x^3y^k\\
  +t^{-8}xy^k.
\end{eqnarray*}

To compute the  term $((+++-+))$ we slide the skein to an area where the
two strands of the twist knot are parallel, as in the first diagram
from Figure~\ref{pppmp1}, then apply Lemma~\ref{strand} at the point 
specified by the arrow, to obtain the first sum in this figure.  

\begin{figure}[h]
  \centering
\scalebox{.30}{\input{pppmp1.pstex_t}}

\caption{}
\label{pppmp1}
\end{figure}
Next, apply Lemma~\ref{strand} at the point specified by the second
arrow in Figure~\ref{pppmp1} to obtain the sum on the second row. 
Resolve the crossings in the first diagram to obtain
that  the first term is  $t^{-4}(y^2-1)(t^2xy^k+A*y^k+xy^{k+1}+t^{-2}xy^k)$. 
The second term is $y\overline{A*y^k}$. Hence 
\begin{eqnarray*}
((+++-+))=(t^{-4}y^2-t^{-4})A*y^k+t^{-6}y\overline{A*y^k}+t^{-4}xy^{k+3}\\+(t^{-2}+t^{-6})xy^{k+2}
-t^{-4}xy^{k+1}+(-t^{-2}-t^{-6})xy^k.
\end{eqnarray*} 

Finally, for $((++++))$, remove the twist and multiply the
skein by $-t^3$,  slide the skein over the top of the diagram
to get the first skein  from Figure~\ref{pppp1} (again only the bottom
of the diagram of the twist knot is shown, and the skein has been moved
to the left off the area where the crossings  occur).
\begin{figure}[h]
  \centering
\scalebox{.30}{\input{pppp1.pstex_t}}

\caption{}
\label{pppp1}
\end{figure}

Apply Lemma~\ref{strand} as specified by first arrow
to obtain the first equality, then apply the lemma again as specified
by second arrow to obtain (after  arranging the terms) the last
sum in Figure~\ref{pppp1}. 
The second term is $-t^{-3}yA*y^k$. After applying Lemma~\ref{strand} as
 specified by the arrow, the first term is
equal to $t^{-1}(-y^2+1)(-t^{-2}yA*y^k-t^{-4}\overline{A*y^k})$. 
So
\begin{eqnarray*}
((++++))=t^{-3}(y^3-2y)A*y^k-t^{-5}(-y^2+1)\overline{A*y^k}. 
\end{eqnarray*}




To simplify the formulas we set $u_i=S_{2i+1}(x)$ and $q=t^4$. 

\begin{theorem}
The action of $(1,-3)_T$ on $K_t({\mathbb T}^2\times I)$ is 
given by
\begin{eqnarray*}
(1,-3)_T\cdot S_k(y)=\sum_{0\leq j\leq 3}t^{2k+2j-1}\alpha_{kj}S_j(y)
\end{eqnarray*}
where
\begin{eqnarray*}
  && \alpha_{0,3}=(-q-1)u_0,\,\alpha_{0,2}=-u_1-2(q+1)u_0,\, \alpha_{0,1}=-2u_1+(q-1)u_0\\&&\alpha_{0,0}=-2u_1+(-6q-1)u_0,\, \alpha_{1,3}=qu_1+(1+q^{-1})u_0,\\
  &&\alpha_{1,2}=(2q+1+q^{-1})u_1+(q+q^{-1}-4q^{-2})u_0,\,\alpha_{1,1}=(q+2)u_1+(q+3)u_0,\\
  &&\alpha_{1,0}=(2q+1+q^{-1})u_1+(2q+4+2q^{-1})u_0, \,\alpha_{2,3}=-qu_2+(q-1)u_1\\
  &&+(-q-2q^{-2})u_0,\, \alpha_{2,2}=(-2q-1)u_2+(-3q-1-q^{-1}-2q^{-2})u_1\\
  &&+(-3q-2-q^{-2}-q^{-3}+q^{-4})u_0,\, \alpha_{2,1}=(-2q-2)u_2+(2q-3-2q^{-1})u_1\\
  &&+(q^3-q^{2}-6q-1-3q^{-1}-3q^{-2}-q^{-3})u_0, \quad \alpha_{2,0}=(-2q-1)u_2\\
  &&+(-4q-5+q^{-1})u_1+(q^{3}-2q^{2}-3q+2-5q^{-1}-2q^{-2})u_0,\\
  &&\alpha_{3,3}=qu_3+u_2+(-1+q^{-1})u_2+(2q^{-3}-3q^{-5})u_0,\quad \alpha_{3,2}=(2q+1)u_3\\
  &&+(q+1+q^{-1})u_2-u_1+(q^{-1}+q^{-2}+q^{-4}-2q^{-5})u_0, \quad \alpha_{3,1}=(2q+2)u_3\\
  &&+(2q+1+2q^{-1})u_2+(-1+2q^{-1}+2q^{-3})u_1+(-q^{-2}+2q^{-3}-q^{-4})u_0,\\
  &&\alpha_{3,0}=(2q+1)u_3+(4q+2+q^{-1})u_2+(6q-1+2q^{-1})u_1\\
  &&+(4q-1-3q^{-1}+q^{-2}-2q^{-3}+q^{-4})u_0.
 \end{eqnarray*}
\end{theorem}

\begin{proof}
  Combining the terms computed above, applying Lemmas~\ref{lemma2} and \ref{lemma6},
  and multiplying with the frame adjusting factor $t^6$,  we obtain:
  \begin{eqnarray*}
    &&(1,-3)_T\cdot y^k= t^{-3}xX_3*y^k+(t^7S_3(y)+t^5S_2(y))A*y^k   \\
    &&\quad +(t^3y+t^5S_2(y))\overline{A*y^k}
    +[t^5S_3(y)+2t^3S_2(y) +(-t^5+2t)S_1(y)\\
      &&\quad +(-2t^3+2t^{-1})]xy^k.
    \end{eqnarray*}
  Substituting the formulas from \S~2 and \S~3, 
switching to the basis $S_j(x)S_k(y)$,  we obtain the formulas from the statement.
  \end{proof}

Let us compute $(1,-2)_T\cdot y^k$, $k=0,1,2,3$. Again we push the $(1,-2)_T$ skein inside
the knot complement and adjust the framing from the plane of the torus to  the plane of the paper, to get the
 skein from Figure~\ref{secondcurve} multiplied by $-t^{9}$. We compute first
the skein from the figure, then adjust  framing. In the figure we  label the 5 crossings
 in the order they are resolved, and   use a boldface
curve for $y^k$, as before.

\begin{figure}[h]
  \centering
\scalebox{.4}{\input{secondcurve.pstex_t}}

\caption{}
\label{secondcurve}
\end{figure}
 As this skein looks similar to  $(1,-3)_T\cdot y^k$,
we expand it in the same way:
\begin{eqnarray*}
(1,-2)_T\cdot y^k=t^4((++++))+t^3((+++-+))+t((+++--))\\
+t((++-))+((+-))+t^{-1}((-)),
\end{eqnarray*}
but now the diagrams of the 6 skeins are different (see Figure~\ref{splusminus}). 
\begin{figure}[h]
\centering
%
\scalebox{.25}{\input{splusminus1.pstex_t}}
\scalebox{.25}{\input{splusminus2.pstex_t}}
\scalebox{.25}{\input{splusminus3.pstex_t}}

\scalebox{.25}{\input{splusminus4.pstex_t}}
\scalebox{.25}{\input{splusminus5.pstex_t}}
\scalebox{.25}{\input{splusminus6.pstex_t}}
 
\caption{}
\label{splusminus}
\end{figure}

To compute $((-))$ we remove the two negative
 twists (and multiply the skein by $t^{-6}$),
 then resolve the two crossings using the Kauffman bracket skein relation.
We obtain
\begin{eqnarray*}
((-))=t^{-6}[t^2x^2y^k+X_3*y^k+X_3*y^k+t^{-2}(-t^2-t^{-2})X_3*y^k]\\
=
t^{-4}x^2y^k+(t^{-6}-t^{-10})X_3*y^k.
\end{eqnarray*}

Next, we focus on $((+-))$. After removing a twist and performing an
isotopy, we obtain the first skein from Figure~\ref{spm}. Now use Lemma~\ref{strand} to obtain the sum in the figure. 
\begin{figure}[h]
  \centering
\scalebox{.25}{\input{spm.pstex_t}}

\caption{}
\label{spm}
\end{figure}

Resolving the crossings with the  skein relation,
we find the first term:
\begin{eqnarray*}
t^{-5}[t^{2}(-t^2-t^{-2})X_4*y^{k+1}+X_4*y^{k+1}+X_4*y^{k+1}+t^{-2}x^2y^{k+1}]\\
=(t^{-5}-t^{-1})X_4*y^{k+1}+t^{-7}xy^{k+1},
\end{eqnarray*}
and the second term:
\begin{eqnarray*}
t^{-7}(t^2x^2y^k+X_3*y^k+X_3*y^k+t^{-2}(-t^2-t^{-2})X_3*y^k)\\
=(t^{-7}-t^{-11})X_3*y^k+t^{-5}x^2y^k.
\end{eqnarray*} 
Combining we obtain
\begin{eqnarray*}
((+-))=(t^{-5}-t^{-1})X_4*y^{k+1}+(t^{-7}-t^{-11})X_3*y^k+t^{-7}xy^{k+1}+t^{-5}x^2y^k.
\end{eqnarray*}

Next we compute $((++-))$, which after an isotopy  becomes the first skein
 in Figure~\ref{sppm}.  
\begin{figure}[h]
  \centering
\scalebox{.25}{\input{sppm.pstex_t}}

\caption{}
\label{sppm}
\end{figure}

Apply Lemma~\ref{strand} to transform this into the sum on
the right side of Figure~\ref{sppm}. Apply Lemma~\ref{strand}
in the first term on the right, then resolve the crossings to obtain that
this term equals 
\begin{eqnarray*}
&&-t^{-2}[(-t^{-2})(t^2x^2y^{k+2}+t^{-2}(-t^2-t^{-2})X_3*y^{k+2}+X_3*y^{k+2}+X_3*y^{k+2}))\\
&&-t^{-4}(t^2(-t^2-t^{-2})X_4*y^{k+1}+t^{-2}x^2y^{k+1}+X_4*y^{k+1}+X_4*y^{k+1})]\\
&&=(t^{-4}-t^{-8})X_3*y^{k+2}+(-t^{-2}+t^{-6})X_4*y^{k+1}+t^{-2}x^2y^{k+2}+t^{-8}x^2y^{k+1}.
\end{eqnarray*}
Then resolve the crossings in the second term to obtain that it is equal
to 
\begin{eqnarray*}
&& -t^{-4}[t^2x^2y^k+t^{-2}(-t^2-t^{-2})X_3*y^k+2X_3*y^k]\\
&&\quad =(-t^{-4}+t^{-8})X_3*y^k-t^{-2}x^2y^k.
\end{eqnarray*}
Combining, we obtain that
\begin{eqnarray*}
  &&((++-))=(t^{-4}-t^{-8})X_3*y^{k+2}+(-t^{-4}+t^{-8})X_3*y^k\\
  &&\quad +(t^{-6}-t^{-2})X_4*y^{k+1} +t^{-2}x^2y^{k+2}+t^{-8}x^2y^{k+1}-t^{-2}x^2 y^k.
\end{eqnarray*}
Which after applying Lemma~\ref{lemma2} for the first and third terms becomes
\begin{eqnarray*}
  ((++-))=(-t^{-4}+t^{-8})X_3*y^k+(t^{-6}-t^{-10})X_2*y^{k+1}\\+t^{-2}x^2y^{k+2}
  +(2t^{-4}-t^{-8})x^2y^{k+1}-t^{-2}x^2y^k.
  \end{eqnarray*}
Compute $((+++--))$ by performing an isotopy over the top of
the knot to obtain the first skein in Figure~\ref{spppmm} (again, ignore
the top of the knot), then  apply Lemma~\ref{strand} as specified by  arrow to obtain the next sum. 
Continue  applying  Lemma~\ref{strand} to each term as specified,
to obtain the sum on the second row in Figure~\ref{spppmm}.  Applying Lemma~\ref{strand} three more times yields

\begin{eqnarray*}
&&((+++--))=-t^{-6}x^2y^{k+2}-t^{-8}x^2y^{k+1}+t^{-6}x^2y^k-t^{-8}x^2y^{k+1}\\
&&\quad \quad -t^{-10}x^2y^k-t^{-10}x^2y^k-t^{-12}Y_1*y^k\\
&&\quad =-t^{-12}Y_1*y^k-t^{-6}x^2y^{k+2}-2t^{-8}x^2y^{k+1}+(t^{-6}-2t^{-10})x^2y^k.
\end{eqnarray*}
\begin{figure}[h]
  \centering
\scalebox{.25}{\input{spppmm.pstex_t}}

\caption{}
\label{spppmm}
\end{figure}

To compute $((+++-+))$, apply Lemma~\ref{strand} as in Figure~\ref{spppmp1}.
The first diagram on the right is $x$ times the second diagram in
Figure~\ref{plusminus} (that denoted by $((+++-+))$ in that figure). 
So the first term on the right is 
\begin{eqnarray*}
(-t^{-6}y^2+t^{-6})xA*y^k-t^{-8}xy\overline{A*y^k}-t^{-4}x^2y^{k+2}-t^{-8}x^2y^{k+2}\\
-t^{-6}x^2y^{k+3}+t^{-4}x^2y^k+t^{-8}x^2y^k+t^{-6}x^2y^{k+1}.
\end{eqnarray*}
\begin{figure}[h]
  \centering
\scalebox{.25}{\input{spppmp1.pstex_t}}

\caption{}
\label{spppmp1}
\end{figure}

The second term on the right can be transformed by an isotopy into the first 
skein in Figure~\ref{spppmp2}. Apply Lemma~\ref{strand} (see arrow)
 to transform it into the sum on the right. Now apply Lemma~\ref{strand}
in each term as specified by the arrows to obtain the sum on the second row. 

\begin{figure}[h]
  \centering
\scalebox{.25}{\input{spppmp2.pstex_t}}
\caption{}
\label{spppmp2}
\end{figure}
The first and  third terms can be combined
into $(-t^{-8}xy-t^{-10}x)Z*y^k$, where $Z*y^k$ is the skein in
Figure~\ref{spppmp3}. We resolve the crossings marked by arrows
using the skein relation to obtain that 
$Z*y^k=t^2xy^{k}+A*y^k+xy^{k+1}+t^{-2}xy^k$. 
The second term from the last sum in Figure~\ref{spppmp2} can be slid back to
the left over the crossing of the twist knot. Then the second term is  just
$-t^{-10} x\overline{A*y^k}$, while the last  is 
$-t^{-12}$ times  the $180^\circ$
rotation of  $C_0*y^k$ in the plane of the paper, so it is in fact
equal to $-t^{-12}C_0*y^k$. Thus
\begin{eqnarray*}
  &&((+++-+))=(-t^{-6}y^2-t^{-8}y-t^{-10}+t^{-6})xA*y^k\\
  &&\quad +(-t^{-8}y-t^{-10})x\overline{A*y^k} -t^{-12}C_0*y^k
  -t^{-6}x^2y^{k+3}\\
  &&\quad +(-t^{-4}-2t^{-8})x^2y^{k+2}-2t^{-10}x^2y^{k+1}+(t^{-4}-t^{-12})x^2y^k.
\end{eqnarray*}

\begin{figure}[h]
  \centering
\scalebox{.25}{\input{spppmp3.pstex_t}}

\caption{}
\label{spppmp3}
\end{figure}

For $((++++))$,   remove the positive twist, multiply by $-t^3$, then slide through the top of the knot diagram
 to obtain the skein in Figure~\ref{spppp}. Apply
Lemma~\ref{strand} as shown by  arrow to  obtain the first
sum, then apply the lemma again in each term, as shown by
arrows, to obtain the second sum. 
\begin{figure}[h]
  \centering
\scalebox{.3}{\input{spppp.pstex_t}}

   \caption{}
  \label{spppp}
\end{figure}
  Apply again Lemma~\ref{strand} in each term to obtain 
\begin{eqnarray*}
&&((++++))=t^{-3}x^2yG*y^k+t^{-5}x^2F*y^k+t^{-5}xyA*y^k\\
  && +t^{-7}x\overline{A*y^k}+t^{-5}x^2G*y^k+t^{-7}xA*y^k +t^{-7}x{A*y^k}+t^{-9}F*y^k\\
  &&=(t^{-3}y+t^{-5})x^2G*y^k+(t^{-5}x^2+t^{-9})F*y^k+(t^{-5}y+2t^{-7})xA*y^k\\
  &&\quad +t^{-7}x\overline{A*y^k}.
\end{eqnarray*}


We set $v_i=S_{2i}(x)$ and $q=t^4$. 

\begin{theorem}
The action of $(1,-2)_T$ on $K_t({\mathbb T}^2\times I)$ is 
given by
\begin{eqnarray*}
(1,-2)\cdot S_k(y)=\sum_{j=0}^3t^{-2k-2j}\beta_{kj} S_j(y),
\end{eqnarray*}
where 
\begin{eqnarray*}
  &&\beta_{0,3}=v_1+qv_0,\, \beta_{0,2}=(q^2+q+1)v_1+(q^2+2q)v_0,\, \beta_{0,1}=(2q^2+2q)v_1\\
  &&+(2q^2+2q)v_0,\,\beta_{0,0}=(2q^2+2q)v_1+(q^2-1-q^{-1})v_0,\\
  &&\beta_{1,3}=-q^{-1}v_2+(q^2-q+1-q^{-1}-q^{-2})v_1+(q^2-1+q^{-1}-q^{-2})v_0,\\
  &&\beta_{1,2}=(-2q^{-1}-q^{-2})v_2+(2q^2-q-q^{-1}-3q^{-2})v_1+(2q^2+2q^{-1})v_0,\\
  &&\beta_{1,1}=(-2q^{-1}-2q^{-2})v_2+(3q^2-3q-5-q^{-1}-q^{-2})v_1\\
  &&+(q^2-3q-7+3q^{-1}-6q^{-2})v_0, \, \beta_{1,0}=(-2q^{-1}-q^{-2})v_2\\
  &&+(-q^2-4q+1-3q^{-1}-4q^{-2})v_1+(-q^2-5q+2-4q^{-1}-4q^{-2})v_0,\\
  &&\beta_{2,3}=q^{-2}v_3+(-q^2+q-q^{-1}+q^{-2}+q^{-3})v_2+(-q^2+1-q^{-1}+q^{-4})v_1\\
  &&+(-1+3q^{-1}-q^{-2}-q^{-3})v_0,\, \beta_{2,2}=(2q^{-1}+q^{-2})v_3+(-2q^3+q^2+q\\
  &&-2-4q^{-1}+2q^{-2}+q^{-3})v_2+(-2q^3+q-4q^{-1}+2q^{-3})v_1+(-q^3-q^2-2\\
  &&+3q^{-1}-q^{-2}-q^{-4})v_0,\\
\end{eqnarray*}
\begin{eqnarray*}
  &&\beta_{2,1}=(2q^{-2}+2q^{-3})v_3 +(-2q^2+2-2q^{-1}+2q^{-2}
  +3q^{-3}+2q^{-4})v_2\\
  &&+(-q^{2}+1+q^{-1}-q^{-2}-q^{-3}+4q^{-4})v_1+(-q^2-q-1+q^{-1}-q^{-2}\\
  &&-q^{-3}+2q^{-4})v_0,\, \beta_{2,0}=(2q^{-2}+q^{-3})v_3+(-2q^2+q+1-2q^{-1}\\
 && +3q^{-2}+3q^{-3}+q^{-4})v_2+(-3q^2-3q+6-q^{-1}-5q^{-2}+3q^{-3}+4q^{-4})v_1\\
  && +(-4q^2-3q+1+3q^{-1}-9q^{-2}-3q^{-3})v_0,\,   \beta_{3,3}=-q^{-3}v_4+(q^{2}-q\\
  &&+q^{-2}-q^{-3}-q^{-4})v_3+(q^2-q+2q^{-2}-4q^{-3}+2q^{-5}+q^{-6}-q^{-7})v_2\\
  &&+(q^{3}+q^2-3q+3-q^{-1}+4q^{-2}-11q^{-3}+5q^{-4}+7q^{-5}-q^{-6}-2q^{-7}-q^{-8})v_1\\ &&(-q^3+3q^2+1+5q^{-1}-q^{-2}-6q^{-3}+2q^{-4}+6q^{-5}-2q^{-6}-2q^{-7}-q^{-8})v_0,\\
  &&\beta_{3,2}=(-2q^{-3}-q^{-4})v_4+(2q^2-q-1+2q^{-2}-2q^{-3}-2q^{-4}-q^{-5})v_3\\
  &&+(-q^3+10q^2-2q-q^{-1}-5q^{-2}-4q^{-3}+2q^{-6}-7q^{-7}-q^{-8})v_2+(-3q^3\\
  &&+24q^2-8q
  +4-3q^{-1}+9q^{-2}-11q^{-3}+q^{-4}-q^{-5}+8q^{_6}-20q^{-7}-q^{-8})v_1\\
  &&+(-2q^3+16q^2-6q  +3-4q^{-1}+6q^{-2}-12q^{-3}-2q^{-5}+6q^{-6}-13q^{-7})v_0,\\
  && \beta_{3,1}=(-2q^{-3}-2q^{-4})v_4+(2q^2-2+2q^{-2}-2q^{-3}-3q^{-4}-2q^{-5})v_3\\
  &&+(-2q^3+3q^2-1-2q^{-1}+5q^{-2}-2q^{-4}+2q^{-6}-2q^{-8})v_2+(-5q^3-2q\\
  && +6+8q^{-1}+6q^{-2}-9q^{-3}-3q^{-4}+6q^{-5}-q^{-7}-5q^{-8})v_1+(-q^2+5\\
  &&+12q^{-1}-5q^{-2}+3q^{-3}-2q^{-4}+4q^{-5})v_0,\, \beta_{3,0}=(-2q^{-3}-q^{-4})v_4+(2q^2-q\\
  &&-1-2q^{-1}+8q^{-2}+10q^{-3}-q^{-4}-2q^{-5}+2q^{-6}-3q^{-7}-q^{-8})v_2+(-3q^3\\
  &&+5q^2+3q+6+9q^{-2}+24q^{-3}+t^{-4}-5t^{-5}+5t^{-6})v_1+(-3q^3+3q^2+3q+6\\
&&  +9q^{-2}+19q^{_3}-4q^{-5}+3q^{-6}+3q^{-7})v_0. 
  \end{eqnarray*}
\end{theorem}

\begin{proof}
  Adding the terms
  and adjusting framing by $-t^9$  we obtain
  \begin{eqnarray*}
    &&(1,-2)_T\cdot y^k=(t^{8}-t^4)X_4*y^{k+1}+(t^6+2t^{-2}-3t^2)X_3*y^k\\
&&    +(-t^4+1)X_2*y^{k+1}
     +t^{-2}Y_1*y^k+[t^{6}y^2+(-t^{8}+t^{4})y\\&&+(-3t^6+t^2)]xA*y^k
      +(t^4y+t^2-t^6)x\overline{A*y^k}+C_0*y^k\\
     &&+(-t^8x^2-t^4)F*y^k
      +(-t^{10}y-t^{8})x^2G*y^k+t^{6}x^2y^{k+3}\\
      &&+3t^4x^2y^{k+2}+(4t^2-2t^6)x^2y^{k+1} +(-3t^4+3)x^2y^k.
    \end{eqnarray*}
  Then use the formulas in \S~2, \S~3
 and switch to the basis $S_j(x)S_k(y)$.
\end{proof}

\end{document}